\documentclass[12pt,twoside]{article}
\usepackage[mathscr]{eucal}
\usepackage{amssymb}
\usepackage{amscd}
\usepackage{amsmath,amsthm}
\usepackage{latexsym}
\usepackage{indentfirst}

\def\C{\mathbb{C}}

\def\I{\mathbb{I}}

\def\ZZ{\mathfrak{Z}}
\def\gg{\mathfrak{g}}

\def\ll{\mathfrak{l}}
\def\pp{\mathfrak{p}}
\def\kk{\mathfrak{k}}

\def\uu{\mathfrak{u}}

\def\oo{\mathfrak{o}}
\def\ss{\mathfrak{s}}

\def\a{\alpha}

\def\d{\delta}
\def\D{\Delta}
\def\e{\varepsilon}
\def\f{\varphi}

\def\pa{\partial}

\def\sm{\setminus}
\def\sub{\subseteq}
\def\ra{\rightarrow}

\def\hra{\hookrightarrow}

\def\lra{\longrightarrow}

\def\cH{{\cal H}}

\def\cD{{\cal D}}

\def\cR{{\cal R}}
\def\cP{{\cal P}}

\def\cO{{\cal O}}

\def\cN{{\cal N}}
\def\cU{{\cal U}}

\def\qqu{\qquad}
\def\qu{\quad}
\def\fa{\forall}

\def\mr{\mathrm}
\def\le{\left}
\def\ri{\right}

\def\tim{\times}
\def\otim{\otimes}
\def\beg{\begin}
\def\ind{\indent}
\def\noind{\noindent}

\def\lb{\linebreak}
\def\ti{\tilde}

\def\en{enumerate}
\def\eq{equation}

\def\hra{\hookrightarrow}

\def\GL{\mr{GL}}
\def\SL{\mr{SL}}
\def\rank{\mr{rank}}

\def\Int{\mr{Int}}

\def\Tr{\mr{Tr}}

\def\End{\mr{End}}

\begin{document}
\title{$K-$structure of $\cU(\gg)$ for $\ss\uu(n,1)$ and $\ss\oo(n,1)$}

\author{Hrvoje Kraljevi\' c, University of Zagreb\thanks{The author was supported by grant no. 4176 of the Croatian Science Foundation and by the QuantiXLie Center of Excellence}}
\date{}
\maketitle

\noind\underline{Abstract.} Let $G$ be the adjoint group of a real simple Lie algebra $\gg_0$ equal\lb either $\ss\uu(n,1)$ or $\ss\oo(n,1),$ $K$ its maximal compact subgroup, $\cU(\gg)$ the universal enveloping algebra of the complexification $\gg$ of $\gg_0$ and $\cU(\gg)^K$ its subalgebra of $K-$invariant elements. By a result of F. Knopp [3] $\cU(\gg)$ is free as a $\cU(\gg)^K-$module, so there exists a $K-$submodule $E$ of $\cU(\gg)$ such that the multiplication defines an isomorphism of $K-$modules $\cU(\gg)^K\otim E\lra\cU(\gg).$ We prove that $E$ is equivalent to the regular representation of $K,$ i.e. that the multiplicity of every $\d\in\hat{K}$ in $E$ equals its dimension. As a consequence we get that for any finitedimensional complex $K-$module $V$ the space $(\cU(\gg)\otim V)^K$ of $K-$invariants is free $\cU(\gg)^K-$module of rank $\dim\,V.$

\section{Introduction}

Let $\gg_0$ be a real simple Lie algebra of noncompact type. Denote by $G$ its adjoint group and choose its maximal compact subgroup $K.$ Let\lb $\gg_0=\kk_0\oplus\pp_0$ be the corresponding Cartan decomposition. Let $\gg,$ $\kk$ and $\pp$ be the complexifications of $\gg_0,$ $\kk_0$ and $\pp_0,$ respectively. Denote by $\cU(\gg)$ and\lb $\cU(\kk)\sub\cU(\gg)$ the universal enveloping algebras of $\gg$ and $\kk.$ Furthermore, denote by $S(\gg)$ and $S(\kk)\sub S(\gg)$ the symmetric algebras over $\gg$ and $\kk$ and by $\cP(\gg)$ and $\cP(\kk)$ the polynomial algebras over $\gg$ and $\kk.$ Then $\cP(\gg)$ and $\cP(\kk)$ can be identified with the symmetric algebras $S(\gg^*)$ and $S(\kk^*)$ over dual spaces $\gg^*$ and $\kk^*$ of $\gg$ and $\kk.$ The Killing form $B$ on $\gg$ allows us to identify $\gg$ with $\gg^*$ and $\kk$ with $\kk^*.$ Thus the algebras $\cP(\gg)$ and $\cP(\kk)$ are identified with $S(\gg)$ and $S(\kk).$ Considering polynomials as complex functions on $\gg$ and $\kk,$ the inclusion $\cP(\kk)\sub\cP(\gg)$ is obtained via the projection $pr:\gg\lra\kk$ along $\pp.$

\newpage

\ind The adjoint action of the group $G$ on $\gg$ extends uniquely to the action by automorphisms on the algebras $\cU(\gg),$ $S(\gg)$ and $\cP(\gg),$ and the subgroup $K$ acts also by automorphisms on the algebras $U(\kk),$ $S(\kk)$ and $\cP(\kk).$ Denote by superscript $G$ (resp. $K)$ the subalgebras of $G-$invariants (resp. $K-$invariants). Then, of course, $ \cU(\gg)^G$ is the center $\ZZ(\gg)$ of $\cU(\gg)$ and $\cU(\kk)^K$ is the center $\ZZ(\kk)$ of $\cU(\kk).$ Obvously, the multiplication defines algebra homomorphisms
$$
\ZZ(\gg)\otim\ZZ(\kk)\ra\cU(\gg)^K,\qu S(\gg)^G\otim S(\kk)^K\ra S(\gg)^K,\qu\cP(\gg)^G\otim\cP(\kk)^K\ra\cP(\gg)^K.
$$
In [3] F. Knopp has proved the following highly nontrivial results:
\beg{tm}\beg{\en}\item[$(a)$] $\ZZ(\gg)\otim\ZZ(\kk)\lra\cU(\gg)^K$ is an isomorphism onto the center of the algebra $\cU(\gg)^K.$
\item[$(b)$] The algebra $\cU(\gg)^K$ is commutative $($i.e. $\cU(\gg)^K=\ZZ(\gg)\ZZ(\kk)\,)$ if and only if $\gg$ is either $\ss\uu(n,1)$ or $\ss\oo(n,1).$ In these cases $\cU(\gg)$ is free as a $\cU(\gg)^K-$module.
\end{\en}
\end{tm}
\ind The symmetrization $\cU(\gg)\lra S(\gg)\simeq\cP(\gg)$ is an isomorphism of vector spaces and of $G-$modules and $(a)$ implies that the homomorphism
$$
\cP(\gg)^G\otim\cP(\kk)^K\lra\cP(\gg)^K
$$
is always injective and by $(b)$ in the cases $\gg=\ss\uu(n,1)$ and $\gg=\ss\oo(n,1)$ this is an isomorphism; furthermore, the last sentence in $(b)$ implies that in these two cases $\cP(\gg)$ is free as a $\cP(\gg)^K-$module.

\section{$K-$harmonic polynomials and the structure of the $\cP(\gg)^K-$module $\cP(\gg)$}

Consider for a while a more general situation. Let $V$ be a complex finite\-dimensional vector space and let $L$ be a closed subgroup of $\GL(V)$ acting fully reducibly on $V.$ Denote by $S(V)$ and $\cP(V)$ the symmetric and the polynomial algebra over $V.$ For $x\in V,$ let $\pa(x):\cP(V)\lra\cP(V)$ be the derivation in the direction $x.$ The map $\pa:V\lra\End(\cP(V))$ extends uniquely to an isomorphism $\pa$ of the symmetric algebra $S(V)$ onto the algebra $\cD(V)$ of linear differential operators on $\cP(V)$ with constant coefficients. Now, one defines the bilinear form $\langle\,\cdot\,,\,\cdot\,\rangle$ on $S(V)\tim\cP(V)$ by
$$
\langle u,f\rangle=[\pa(u)f](0),\qqu u\in S(V),\,\,f\in\cP(V).
$$
This is a pairing, i.e. nondegenerate in each variable. Now, consider the sub\-algebras of $L-$invariants $S(V)^L$ and $\cP(V)^L$ and their maximal ideals (of codimension 1)
$$
S_+(V)^L=\bigoplus_{k>0}S^k(V)^L,\qu\cP_+(V)^L=\bigoplus_{k>0}\cP^k(V)^L=\{f\in\cP(V)^L;\,\,f(0)=0\}.
$$
Define the (graded) space of so called $L-$harmonic polynomials on $V:$
$$
\cH_L(V)=\{f\in\cP(V);\,\,\pa(u)f=0\,\,\,\fa u\in S_+(V)^L\}.
$$
As noticed in [4] and [5] the obvious equality
$$
\langle uv,f\rangle =\langle u,\pa(v)f\rangle,\qqu u,v\in S(V),\qu f\in\cP(V),
$$
implies easily that
$$
\cH_L(V)=\{f\in\cP(V);\,\,\langle u,f\rangle=0\,\,\,\fa u\in S(V)S_+(V)^L\}.
$$
Part of the Helgason's results in [2] (see also Propositions 3 and 4 in [4]) can be stated as follows:
\beg{p} Suppose that the group $L$ is connected and that there exists an $L-$invariant symmetric bilinear form $B:V\tim V\lra\C$ and a real form $V_0$ of $V$ such that the restriction of $B$ to $V_0\tim V_0$ is a scalar product and that the group $L$ is the complexification of its subgroup $L_0=\{g\in L;\,\,gV_0=V_0\}.$ Then
$$
\cP(V)=\cP(V)\cP_+(V)^L\oplus\cH_L(V).
$$
\end{p}
\ind Note that the conditions on the pair $(L,V)$ in Proposition 1 are obviously satisfied for the action of the complexification $K^{\C}$ of the group $K$ on $\gg,$\lb especially in the cases $\gg_0=\ss\uu(n,1)$ and $\gg_0=\ss\oo(n,1).$\\
\ind Consider any subgroup $L\sub\GL(V)$ acting fully reducibly on a finitedimensional complex vector space $V.$ If $N$ is any graded subspace of $\cP(V)$ such that
\beg{\eq}
\cP(V)=\cP(V)\cP_+(V)^L\oplus N
\end{\eq}
then it is easy to see (Proposition 1 in [4]) that the multiplication defines a surjective map
\beg{\eq}
\cP(V)^L\otim N\lra\cP(V).
\end{\eq}
Kostant's Lemma 1 in [4] can be stated as follows:
\beg{p} The following properties are mutally equivalent:
\beg{\en}
\item[$(a)$] For every $N,$ such that $(1)$ holds true, the map $(2)$ is also injective, i.e. an isomorphism.
\item[$(b)$] For some $N,$ such that $(1)$ holds true, the map $(2)$ is injective.
\item[$(c)$] $\cP(V)$ is free as a $\cP(V)^L-$module.
\end{\en}
\end{p}
\ind Thus, by the last sentence in $(b)$ of Theorem 1 we get from Propositions 1 and 2:
\beg{tm} For $\gg=\ss\uu(n,1)$ and for $\gg=\ss\oo(n,1)$ we have:
\beg{\en}
\item[$(a)$] $\cP(\gg)=\cP(\gg)\cP_+(\gg)^K\oplus\cH_K(\gg).$
\item[$(b)$] The multiplication defines an isomorphism $\cP(\gg)^K\otim\cH_K(\gg)\simeq\cP(\gg).$
\end{\en}
\end{tm}

\section{The $K-$module of $K-$harmonic polynomials}

\ind Let $\cN$ be the zero set in $\gg$ of the ideal $\cP(\gg)\cP_+(\gg)^K$ generated by $\cP_+(\gg)^K$ in $\cP(\gg):$
$$
\cN=\{x\in\gg;\,\,f(x)=0\,\,\,\fa f\in\cP(\gg)\cP_+(\gg)^K\}=\{x\in\gg;\,\,f(x)=0\,\,\,\fa f\in\cP_+(\gg)^K\}.
$$
\ind By Proposition 16 in [4] the zero set
$$
\cN_G=\{x\in\gg;\,\,f(x)=0\,\,\,\fa f\in\cP_+(\gg)^G\}
$$
is exactly the set of all nilpotent elements in the Lie algebra $\gg.$ Analogously
$$
\cN_K=\{x\in\kk;\,\,f(x)=0\,\,\,\fa f\in\cP_+(\kk)^K\}
$$
is the set of all nilpotent elements in the reductive Lie algebra $\kk.$ Now, $\cP(\gg)^K=\cP(\gg)^G\otim\cP(\kk)^K$ by the Knopp's theorem, so we get
\beg{p} $\cN$ is the set of all nilpotent elements in $\gg$ whose projection to $\kk$ along $\pp$ is nilpotent in the reductive Lie algebra $\kk:$
$$
\cN=\{x\in\gg;\,\,x\in\cN_G,\,\,\,pr\,x\in\cN_K\}.
$$
\end{p}
\ind We call the elements of $\cN$ $K-$nilpotent elements in $\gg.$\\

\ind By the Harish$-$Chandra isomorphism and by the Chevalley's theorem on Weyl group invariants we know that the algebra $\cP(\gg)^G$ is generated by $\ell=\rank\,\gg$ homogeneous algebraically independent $G-$invariant polynomials $f_1,\ldots,f_{\ell}$ and the algebra $\cP(\kk)^K$ is generated by $k=\rank\,\kk$ homogeneous algebraically independent $K-$invariant polynomials $\f_1,\ldots,\f_k.$ Since in the cases $\gg_0=\ss\uu(n,1)$ and $\gg_0=\ss\oo(n,1)$
$$
\cP(\gg)^K=\cP(\gg)^G\cP(\kk)^K\simeq\cP(\gg)^G\otim\cP(\kk)^K,
$$
the algebra $\cP(\gg)^K$ is generated by $\ell+k$ homogeneous algebraically independent polynomials $f_1,\ldots,f_{\ell},\f_1,\ldots,\f_k.$ Thus,
$$
\cN=\{x\in\gg;\,\,f_1(x)=\cdots=f_{\ell}(x)=\f_1(x)=\cdots=\f_k(x)=0\},
$$
so the set $\cN$ is a Zariski closed subset of $\gg$ of dimension
$$
\dim\,\cN=\dim\,\gg-\ell-k.
$$
More generally, for any $(\xi,\eta)=(\xi_1,\ldots,\xi_{\ell},\eta_1,\ldots,\eta_k)\in\C^{\ell+k}$ we define a $K^{\C}-$stable Zariski closed subset $\cN(\xi,\eta)$ of $\gg:$
$$
\cN(\xi,\eta)=\{x\in\gg;\,\,f_j(x)=\xi_j,\,j=1,\ldots,\ell,\,\,\f_i(x)=\eta_i,\,i=1,\ldots,k\}.
$$
Obviously,
$$
\dim\,\cN(\xi,\eta)=\dim\,\gg-\ell-k,\qqu(\xi,\eta)\in\C^{\ell+k}.
$$
As in [4] and [5] we conclude from Theorem $2(a):$
\beg{p} The restriction of polynomials in $\cP(\gg)$ to the set $\cN(\xi,\eta)$ induces an isomorphism of $K-$modules
$$
\cH_K(\gg)\simeq\cP(\cN(\xi,\eta))=\cR(\cN(\xi,\eta)),\qqu(\xi,\eta)\in\C^{k+\ell}.
$$
\end{p}
Here for any subset $S\sub\gg$ we set
$$
\cP(S)=\{f|S;\,\,f\in\cP(\gg)\}
$$
and for any algebraic variety $S$ $\cR(S)$ denotes the algebra of regular functions on $S.$\\
\ind The dimensions and the ranks $\ell=\rank\,\gg$ and $k=\rank\,\kk$ in our cases are the following:
$$
\beg{array}{|c|c|c|c|c|}
\hline
\gg&\dim\gg&\dim\kk&\ell&k\\
\hline
\ss\uu(n,1)&n^2+2n&n^2&n&n\\
\hline
\ss\oo(2n,1)&2n^2+n&2n^2-n&n&n\\
\hline
\ss\oo(2n+1,1)&2n^2+3n+1&2n^2+n&n+1&n\\
\hline
\end{array}
$$
So we see that in each case
\beg{\eq}
\dim\,\cN(\xi,\eta)=\dim\,\kk=\dim\,K^{\C},\qqu(\xi,\eta)\in\C^{\ell+k}.
\end{\eq}

\noind{\bf Remark:} By the exercise 13) in \S 13 in [1] (p. 268) we can choose the following generators $f_i,\f_j$ od $\cP(\gg)^K:$
\beg{\en}\item[$(a)$] For $\gg_0=\ss\uu(n,1)$
$$
f_i(x)=\Tr\,x^{i+1},\qu1\leq i\leq n,\qqu\f_j(x)=\Tr\,(pr\,x)^j,\qu1\leq j\leq n.
$$
\item[$(b)$] For $\gg_0=\ss\oo(2n,1)$
$$
f_i(x)=\Tr\,x^{2i},\qu1\leq i\leq n,\qqu\f_j(x)=\Tr\,(pr\,x)^{2j},\qu1\leq j\leq n-1,
$$
$$
\f_n(x)^2=(-1)^n\det\,(pr\,x).
$$
\item[$(c)$] For $\gg_0=\ss\oo(2n+1,1)$
$$
f_i(x)=\Tr\,x^{2i},\qu1\leq i\leq n,\qqu f_{n+1}(x)^2=(-1)^{n+1}\det\,x,
$$
$$
\f_j(x)=\Tr\,(pr\,x)^{2j},\qu1\leq j\leq n.
$$
\end{\en}

\ind Consider the action of the complex group $K^{\C}$ on $\gg.$ For $x\in\gg$ denote by $\cO_x$ its $K^{\C}-$orbit. Then of course
\beg{\eq}
\dim\,\cO_x=\dim\,K^{\C}/K_x^{\C}=\dim\,K^{\C}-\dim\,K_x^{\C},
\end{\eq}
where $K_x^{\C}$ denotes the stabilizer of the point $x$ in the group $K^{\C}.$ So, if $K_x^{\C}$ is trivial
\beg{\eq}
\dim\,\cO_x=\dim\,K^{\C}=\dim\,\cN(\xi,\eta).
\end{\eq}
\beg{lm} There exists $x\in\gg$ such that the stabilizer $K_x^{\C}$ is trivial. In this case let $(\xi,\eta)=(f_1(x),\ldots,f_{\ell}(x),\f_1(x),\ldots,\f_k(x)).$ The orbit $\cO_x$ is open in $\cN(\xi,\eta).$
\end{lm}

\ind We prove this Lemma in Section 4.\\

\ind Let $x\in\gg$ be as in Lemma 1, i.e. such that its stabilizer in $K^{\C}$ is trivial. Set
$$
(\xi,\eta)=(f_1(x),\ldots,f_{\ell}(x),\f_1(x),\ldots,\f_k(x))\in\C^{\ell+k}.
$$
We know that $\dim\,\cO_x=\dim\,\cN(\xi,\eta),$ so the $K^{\C}-$orbit $\cO_x$ is open in $\cN(\xi,\eta).$ Thus, the restriction to $\cO_x$ is an isomorphism of $\cP(\cN(\xi,\eta))=\cR(\cN(\xi,\eta))$ onto $\cP(\cO_x).$ Now, if the algebraic variety $\cN(\xi,\eta)$ would be irreducible and if we would have
\beg{\eq}
\dim\,\cN(\xi,\eta)\sm\cO_x\leq\dim\cN(\xi,\eta)-2,
\end{\eq}
(this holds true in the settings of [4] and [5] since the dimensions of all the orbits have the same parity) we could conclude by a theorem from algebraic geometry that $\cP(\cO_x)=\cR(\cO_x)\simeq\cR(K^{\C})$ as $K^{\C}-$modules and by the Frobenius reciprocity we could get that the multiplicity $m(\d)$ of any irreducible finitedimensional representation $\d$ of $K^{\C}$ in the $K^{\C}-$module $\cH_K(\gg)\simeq\cR(\cO_x)$ equals its dimension $d(\d).$ Unfortunately, $(6)$ is not true. In fact, in the case $\gg=\ss\uu(n,1)$ the algebraic set $\cN=\cN(0,0)$ is even not irreducible $-$ there exist two open orbits in $\cN,$ and in the complement of these two orbits there exist orbits of dimension $\dim\,\cN-1.$ In the case $\gg=\ss\oo(n,1),$ $n\geq3,$ there also exist $K^{\C}-$orbits in $\cN(\xi,\eta)$ of dimension $\cN(\xi,\eta)-1.$\\
\ind So, we get only the inclusion of $K-$modules $\cH_K(\gg)\hra\cR(K^{\C})$ and we may conclude only that
\beg{\eq}
m(\d)\leq d(\d)
\end{\eq}
for every irreducible finitedimensional representation $\d$ of $K.$ In fact, the equality holds true although we do not know {\it a priory} that $\cP(\cO_x)=\cR(\cO_x);$ it comes out {\it a posteriory}:
\beg{tm} The multiplicity of every irreducible finitedimensional representation $\d$ of the compact group $K$ in the $K-$module $\cH_K(\gg)$ of $K-$harmonic polynomials on $\gg$ is equal to its dimension $d(\d).$
\end{tm}

\ind To prove Theorem 3 we use the compact form $K$ of the complex group $K^{\C}.$ Denote by $\cP(Kx)$ the restriction of the polynomial algebra $\cP(\gg)$ to the $K-$orbit $Kx.$ Note that the fact that $K^{\C}$ is the complexification of $K$ easily implies that the restriction $\cO_x\ra Kx$ induces an isomorphism of $K-$modules $\cP(\cO_x)$ onto $\cP(Kx).$ Thus, as a $K-$module we have
\beg{\eq}
\cP(Kx)=\bigoplus_{\d\in\hat{K}}m(\d)\d.
\end{\eq}
The subalgebra $\cP(Kx)$ of the algebra $C(Kx)$ of all complex continuous functions on the compact space $Kx$ evidently distinguishes the points of $Kx.$ Furthermore, this subalgebra is closed under complex conjugation. This is implied by the fact that the set $Kx$ is contained in a real form of the complex vector space $\gg.$ This follows from the fact that the compact group K is contained in a maximal compact subgroup U of the complex group $G^{\C}=\Int(\gg)$ and the Lie algebra $\uu$ of $U$ is a real form of $\gg.$ Finally, the algebra $\cP(Kx)$ obviously contains constants. Thus, by the Stone$-$Weierstrass theorem the subalgebra $\cP(Kx)$ is uniformly dense in $C(Kx).$ Now, the Peter$-$Weyl theo\-rem implies that $m(\d)=d(\d)$ for all $\d\in\hat{K}.$ This proves Theorem 3.\\

\ind The symmetrization $U(\gg)\lra S(\gg)\simeq\cP(\gg)$ is an $K-$module isomor\-phism. Let $H_K$ be the inverse image of $\cH_K(\gg)$ in $U(\gg).$ The immediate consequence of Theorems 2 and 3 is
\beg{tm} The multiplication induces an isomorphism of $K-$modules $U(\gg)^K\otim H_K\simeq U(\gg).$ The multiplity of every $\d\in\hat{K}$ in the $K-$module $H_K$ is equal to its dimension $d(\d).$
\end{tm}
\beg{ko} Let $V$ be a finitedimensional $K-$module. Then the space of $K-$invariants $(U(\gg)\otim V)^K$ is a free $U(\gg)^K-$module of finite rank $\dim\,V.$
\end{ko}
\ind By Theorem 4 we have
$$
(U(\gg)\otim V)^K\simeq(U(\gg)^K\otim H_K\otim V)^K=U(\gg)^K\otim(H_K\otim V)^K.
$$
Thus, $U(\gg)$ is a free $U(\gg)^K-$module of rank $\dim\,(H_K\otim V)^K.$ Now, let $n(\e)$ be the multiplicity of $\e\in\hat{K}$ in $V.$ Then
$$
(H_K\otim V)^K\simeq\le(\le(\oplus_{\d\in\hat{K}}d(\d)\d\ri)\otim\le(\oplus_{\e\in\hat{K}}n(\e)\e\ri)\ri)^K=\oplus_{\d,\e\in\hat{K}}d(\d)n(\e)(\d\otim\e)^K,
$$
so
$$
\dim\,(H_K\otim V)^K=\sum_{\d,\e\in\hat{K}}d(\d)n(\e)\dim\,(\d\otim\e)^K.
$$
By the Schur's lemma $\dim\,(\d\otim\e)^K$ is $1$ if $\d$ and $\e$ are contragredient to each other and $0$ otherwise. Since the dimensions of contragredient representations are equal, we get
$$
\dim\,(H_K\otim V)^K=\sum_{\d\in\hat{K}}n(\d)d(\d)=\dim\,V.
$$

\section{Proof of Lemma 1}

(1) $\gg_0=\ss\uu(n,1).$ We realize this Lie algebra as
$$
\gg_0=\{A\in\ss\ll(n+1,\C);\,\,A^*=-\Gamma A\Gamma\},
$$
where $\Gamma=\mr{diag}\,(1,\ldots,1,-1).$ Then $\gg=\ss\ll(n+1,\C)$ and $K^{\C}=\ti{K}^{\C}/Z,$ where $Z=\{\mr{diag}\,(\a,\ldots,\a);\,\,\a^{n+1}=1\}$ is the center of $\SL(n+1,\C)$ and
$$
\ti{K}^{\C}=\le\{\le[\beg{array}{cc}B&0\\0&(\det\,B)^{-1}\end{array}\ri];\,\,B\in\GL(n,\C)\ri\}.
$$
Now, we can take for $x$ the elementary $(n+1)\tim(n+1)$ Jordan block:
$$
x=\le[\beg{array}{ccccc}0&1&\cdots&0&0\\0&0&\cdots&0&0\\\vdots&\vdots&&\vdots&\vdots\\0&0&\cdots&0&1\\0&0&\cdots&0&0\end{array}\ri].
$$
The centralizer $M_x$ of $x$ in the algebra of all $(n+1)\tim(n+1)$ matrices consists of  all polynomials in $x,$ i.e.
$$
M_x=\le\{\le[\beg{array}{cccccc}\a_0&\a_1&\a_2&\cdots&\a_{n-1}&\a_n\\0&\a_0&\a_1&\cdots&\a_{n-2}&a_{n-1}\\0&0&\a_0&\cdots&\a_{n-3}&\a_{n-2}\\\vdots&\vdots&\vdots&&\vdots&\vdots\\0&0&0&\cdots&\a_0&\a_1\\0&0&0&\cdots&0&\a_0\end{array}\ri];\,\,\a_0,\a_1,\ldots,\a_n\in\C\ri\}.
$$
So, we conclude that the centralizer of $x$ in $\ti{K}^{\C}$ is precisely the center $Z$ of $\SL(n+1,\C),$ thus the stabilizer of $x$ in $K^{\C}$ is trivial.\\
\ind(2) $\gg_0=\ss\oo(2n+1,1).$ We choose the following realizations:
$$
\gg=\ss\oo(2n+2,\C)=\{A\in\gg\ll(2n+2,\C);\,\,A^t=-\Gamma A\Gamma\},
$$
$$
\kk=\le\{\le[\beg{array}{cc}B&0\\0&0\end{array}\ri];\,\,B\in\gg\ll(2n+1,\C),\,\,B^t=-\Gamma_0B\Gamma_0\ri\}.
$$
Here the superscript $t$ denotes the matrix transpose and
$$
\Gamma_0=\le[\beg{array}{ccc}0&I_n&0\\I_n&0&0\\0&0&1\end{array}\ri],\qqu\Gamma=\le[\beg{array}{cc}\Gamma_0&0\\0&1\end{array}\ri],
$$
$I_n$ being the $n$ by $n$ identity matrix. Denoting as usual the space of all $n\tim m$ complex matrices by $M_{n,m}(\C)$ and $M_n(\C)=M_{n,n}(\C),$ we have
$$
\kk=\le\{\le[\beg{array}{cccc}A&B&a&0\\C&-A^t&b&0\\-b^t&-a^t&0&0\\0&0&0&0\end{array}\ri];\,\,A,B,C\in M_n(\C),\,\,B^t=-B,\,\,C^t=-C,\,\,a,b\in M_{n,1}(\C)\ri\}
$$
and
$$
\gg=\le\{X+\le[\beg{array}{cccc}0&0&0&c\\0&0&0&d\\0&0&0&\a\\-d^t&-c^t&-\a&0\end{array}\ri];\,\,X\in\kk,\,\,c,d\in M_{n,1}(\C),\,\,\a\in\C\ri\}.
$$
Furthermore,
$$
K^{\C}=\le\{\le[\beg{array}{cc}A&0\\0&1\end{array}\ri];\,\,A\in\SL(2n+1,\C),\,\,A^t\Gamma_0A=\Gamma_0\ri\}.
$$
Let $J$ denote the elementary $n$ by $n$ Jordan block and let $e_j\in M_{n,1}(\C)$ be the column matrix with $1$ in the $j-$th row and zeros elsewhere. Set
$$
x=\le[\beg{array}{cccc}J&0&e_n&0\\0&-J^t&0&e_1\\0&-e_n^t&0&0\\-e_1^t&0&0&0\end{array}\ri].
$$
This element of $\gg$ is an invertible matrix which is up to the change of some signs (to be precise, on the places $1,n+1,n+2,\ldots,2n)$ the matrix of the following cyclic permutation of the set of indices $\{1,2,\ldots,2n+2\}:$
$$
1\ra2\ra\cdots\ra n-1\ra n\ra2n+1\ra2n+1\ra n+1\ra n+2\ra\cdots\ra2n\ra1.
$$
Thus, we conclude that the stabilizer (i.e. the centralizer) of $x$ in $K^{\C}$ is trivial.

\newpage

\ind (3) $\gg=\ss\oo(2n,1).$ We choose the following realizations
$$
\gg=\ss\oo(2n+1,\C)=\{A\in\gg\ll(2n+1,\C);\,\,A^{t}=-\Gamma A\Gamma\},
$$
$$
\kk=\le\{\le[\beg{array}{cc}B&0\\0&0\end{array}\ri];\,\,B\in\gg\ll(2n,\C),\,\,B^{t}=-\Gamma_0B\Gamma_0\ri\},
$$
$$
\Gamma_0=\le[\beg{array}{cc}0&I_n\\I_n&0\end{array}\ri],\qqu\Gamma=\le[\beg{array}{cc}\Gamma_0&0\\0&1\end{array}\ri],
$$
Then
$$
\kk=\le\{\le[\beg{array}{ccc}A&B&0\\C&-A^t&0\\0&0&0\end{array}\ri];\,\,A,B,C\in M_n(\C),\,\,B^t=-B,\,\,C^t=-C\ri\}
$$
and
$$
\gg=\le\{X+\le[\beg{array}{ccc}0&0&a\\0&0&b\\-b^t&-a^t&0\end{array}\ri];\,\,X\in\kk,\,\,a,b\in M_{n,1}(\C)\ri\}.
$$
As in $(2)$ let $J$ denote the elementary $n$ by $n$ Jordan block and let 
$$
\D=\le[\beg{array}{ccccc}0&\cdots&0&0&0\\\vdots&&\vdots&\vdots&\vdots\\0&\cdots&0&0&0\\0&\cdots&0&0&1\\0&\cdots&0&-1&0\end{array}\ri]\in M_n(\C).
$$
The matrix
$$
x_0=\le[\beg{array}{ccc}J&\D&0\\0&-J^t&0\\0&0&0\end{array}\ri],
$$
is a representative of the $K^{\C}-$orbit of all principal nilpotent elements of $\kk.$ By the Kostant's results in [3] the stabilizer $K_{x_0}^{\C}$ of $x_0$ in $K^{\C}$ is an $n-$dimensional connected simply connected unipotent subgroup whose Lie algebra is the centralizer $\kk_{x_0}$ of $x_0$ in $\kk.$\\
\ind $(3a)$ Suppose first that $n$ is odd, $n=2k+1.$ By solving a system of linear equations one finds that $\kk_{x_0}$ consists of all matrices of the form
\beg{\eq}
\le[\beg{array}{ccc}A&B&0\\0&-A^t&0\\0&0&0\end{array}\ri],
\end{\eq}
where $B$ is $n$ by $n$ antisymmetric matrix such that for some $\a_1,\a_2,\ldots,\a_n\in\C$ its first row is
$$
\le[\beg{array}{ccccccccc}0&\a_1&0&\a_2&0&\cdots&0&\a_k&\a_{k+1}\end{array}\ri],
$$
its last column is
$$
\le[\beg{array}{ccccccccc}\a_{k+1}&\a_{k+2}&0&\a_{k+3}&0\cdots&0&\a_{2k+1}&0\end{array}\ri]^t,
$$
the inner entries of $B$ are either $0,$ or $\pm\a_j,$ $2\leq j\leq k,$ or $\pm2\a_j,$ $k+2\leq j\leq 2k,$ and $A$ is a strictly upper triangular $n$ by $n$ matrix whose first row is
$$
\le[\beg{array}{ccccccccc}0&\a_{2k+1}&0&\a_{2k}&0&\cdots&0&\a_{k+2}&-\a_{k+1}\end{array}\ri]
$$
and every paralel with the main diagonal is constant (i.e. $A$ is a polynomial in $J).$ E.g. for $n=7$ $(k=3)$
$$
A=\le[\beg{array}{ccccccc}0&\a_7&0&\a_6&0&\a_5&-\a_4\\0&0&\a_7&0&\a_6&0&\a_5\\0&0&0&\a_7&0&\a_6&0\\0&0&0&0&\a_7&0&\a_6\\0&0&0&0&0&\a_7&0\\0&0&0&0&0&0&\a_7\\0&0&0&0&0&0&0\end{array}\ri],
$$
$$ B=\le[\beg{array}{ccccccc}0&\a_1&0&\a_2&0&\a_3&\a_4\\-\a_1&0&-\a_2&0&-\a_3&0&\a_5\\0&\a_2&0&\a_3&0&-2\a_5&0\\-\a_2&0&-\a_3&0&2\a_5&0&\a_6\\0&\a_3&0&-2\a_5&0&-2\a_6&0\\-\a_3&0&2\a_5&0&2\a_6&0&\a_7\\-\a_4&-\a_5&0&-\a_6&0&-\a_7&0\end{array}\ri].
$$
\ind $(3b)$ Consider now the case of $n$ even, $n=2k.$ As in $(3a)$ one finds that $\kk_{x_0}$ consists of all matrices of the form $(9)$ where $B$ is $n$ by $n$ antisymmetric matrix whose first row is
$$
\le[\beg{array}{cccccccc}0&\a_1&0&\a_2&0&\cdots&0&\a_k\end{array}\ri],
$$
its last column is
$$
\le[\beg{array}{cccccccccc}\a_k&0&\a_{k+2}&0&\a_{k+3}&0&\cdots&0&\a_{2k}&0\end{array}\ri]^t,
$$
the inner entries of its antidiagonal are $\pm\a_{k+1},$ all the other inner entries are either $0,$ or $\pm\a_j,$ $2\leq j\leq k-1,$ or $\pm2\a_j,$ $k+2\leq j\leq2k-1,$ and $A$ is the strictly upper triangular $n$ by $n$ matrix whose first row is
$$
\le[\beg{array}{cccccccccc}0&\a_{2k}&0&\a_{2k-1}&0&\cdots&0&\a_{k+2}&0&\a_{k+1}-\a_k\end{array}\ri]
$$
and every paralel with the main diagonal is constant. E.g. for $n=6$ $(k=3)$
$$
A=\le[\beg{array}{cccccc}0&\a_6&0&\a_5&0&\a_4-\a_3\\0&0&\a_6&0&\a_5&0\\0&0&0&\a_6&0&\a_5\\0&0&0&0&\a_6&0\\0&0&0&0&0&\a_6\\0&0&0&0&0&0\end{array}\ri],
$$
$$
B=\le[\beg{array}{cccccc}0&\a_1&0&\a_2&0&\a_3\\-\a_1&0&-\a_2&0&-\a_4&0\\0&\a_2&0&\a_4&0&\a_5\\-\a_2&0&-\a_4&0&-2\a_5&0\\0&\a_4&0&2\a_5&0&\a_6\\-\a_3&0&-\a_5&0&-\a_6&0\end{array}\ri].
$$
\ind Now, since $\pp$ is $K^{\C}-$stable, for any $y\in\pp$ the stabilizer (resp. the centra\-lizer) of $x=x_0+y$ in $K^{\C}$ (resp. $\kk)$ is the stabilizer (resp. the centralizer) of $y$ in $K_{x_0}^{\C}$ (resp. $\kk_{x_0}).$ Let us compute the centralizer of
$$
y=\le[\beg{array}{ccc}0&0&0\\0&0&e_1\\-e_1^t&0&0\end{array}\ri]\in\pp^{\C}
$$
in $\kk_{x_0}.$ An element $(9)$ of $\kk_{x_0}$ centralizes $y$ if and only if
$$
Be_1=0\qqu\mr{and}\qqu A^te_1=0.
$$
Now, in the case $(3a)$ we have
$$
Be_1=\le[\beg{array}{cccccccc}0&-\a_1&0&-\a_2&0&\cdots&-\a_k&-\a_{k+1}\end{array}\ri]^t,
$$
$$
A^te_1=\le[\beg{array}{ccccccccc}0&\a_{2k+1}&0&\a_{2k}&0&\cdots&0&\a_{k+2}&-\a_{k+1}\end{array}\ri]^t,
$$
and in the case $(3b)$
$$
Be_1=\le[\beg{array}{cccccccc}0&-\a_1&0&-\a_2&0&\cdots&0&-\a_k\end{array}\ri]^t,
$$
$$
A^te_1=\le[\beg{array}{cccccccccc}0&\a_{2k}&0&\a_{2k-1}&0&\vdots&0&\a_{k+2}&0&\a_{k+1}-\a_k\end{array}\ri]^t.
$$
In both cases we conclude that $(5)$ is in the centralizer of $y$ in $\kk_{x_0}$ if and only if $\a_1=\a_2=\cdots=\a_n=0,$ i.e. if and only if $A=B=0.$ Thus,
$$
x=x_0+y=\le[\beg{array}{ccc}J&\D&0\\0&-J^t&e_1\\-e_1^t&0&0\end{array}\ri]
$$
is an element of $\gg$ whose stabilizer in $K^{\C}$ is trivial. This completes the proof of Lemma 1.

\beg{thebibliography}{5}
\bibitem{pa} N. Bourbaki, \emph{Lie Groups and Lie Algebras, Chapters 7$-$9}, Springer$-$Verlag, Berlin$-$Heidelberg, 2005.
\bibitem{pa} S. Helgason, \emph{Some results in invariant theory}, Bulletin of the American Mathematical Society, vol. 68 (1962), pp 367$-$371.
\bibitem{pa} F. Knopp, \emph{Der Zentralisator einer Liealgebra in einer einh\" ullenden Algebra}, Journal f\" ur die reine und angewandte Mathematik, vol. 406 (1990), pp 5$-$9.
\bibitem{pa} B. Kostant, \emph{Lie group representations on polynomial rings}, American Journal of Mathematics, vol. 86 (1963), pp 327$-$402.
\bibitem{pa} B. Kostant and S. Rallis, \emph{Orbits and representations associated with symmetric spaces}, American Journal of Mathematics, vol. 93 (1971), pp 753$-$809.

\end{thebibliography}

\end{document}